\documentclass[journal]{IEEEtran}
\usepackage[T1]{fontenc}
\usepackage{amsmath,amssymb,amsthm,booktabs,graphicx,microtype}
\usepackage{cite}
\usepackage{array}
\usepackage{algorithm}
\usepackage[noend]{algpseudocode}
\usepackage[hidelinks]{hyperref}
\usepackage{xurl}
\graphicspath{{Figures/}}
\hypersetup{pdftitle={Joint Reservation of Electric Vehicle Charging Ports: Who Gains and When to Batch},pdfauthor={Luca Ambrosino, Hong T. M. Chu, Laurent El Ghaoui and Giuseppe C. Calafiore}}
\newtheorem{theorem}{Theorem}
\newtheorem{proposition}{Proposition}
\newcommand{\WashingtonSaving}{4.63}
\newcommand{\ChicagoSaving}{8.92}
\newcommand{\ManhattanSaving}{2.03}
\newcommand{\HighGreedyService}{93.84}
\newcommand{\HighGap}{6.16}
\newcommand{\HighGapLo}{5.70}
\newcommand{\HighGapHi}{6.68}
\newcommand{\HighSaving}{3.95}
\newcommand{\HighUnservedEnergy}{2.47}
\newcommand{\FallbackFifty}{0.29}
\newcommand{\FallbackHundred}{3.37}
\newcommand{\FallbackTwoHundred}{9.53}
\begin{document}
\title{Joint Reservation of Electric Vehicle Charging Ports: Who Gains and When to Batch}
\author{Luca~Ambrosino, Hong~T.~M.~Chu, Laurent~El~Ghaoui, and Giuseppe~C.~Calafiore%
\thanks{This research did not receive any specific grant from funding agencies in the public, commercial, or not-for-profit sectors.}%
\thanks{L. Ambrosino is with the Department of Electronics and Telecommunications, Politecnico di Torino, 10129 Torino, Italy (corresponding author, e-mail: \texttt{luca.ambrosino@polito.it}).}%
\thanks{H. T. M. Chu and L. El Ghaoui are with the College of Engineering and Computer Science and the Center for Environmental Intelligence, VinUniversity, Hanoi, Vietnam (e-mail: \texttt{hong.ctm@vinuni.edu.vn}; \texttt{laurent.eg@vinuni.edu.vn}).}%
\thanks{G. C. Calafiore is with the Department of Electronics and Telecommunications, Politecnico di Torino, 10129 Torino, Italy, and with the Center for Environmental Intelligence, VinUniversity, Hanoi, Vietnam (e-mail: \texttt{giuseppe.calafiore@polito.it}).}}
\maketitle

\begin{abstract}
We study how a charging platform should reserve ports for a batch of electric vehicles whose drivers have individual travel limits, connector compatibility and charging budgets, serving as many as possible at the least combined cost of travel, electricity, waiting and charging time. We show that, because each vehicle requests a fixed amount of energy from a single port, the budget constraints amount to screening out unaffordable vehicle-port pairs, so the integer program is exactly a rectangular assignment problem solvable in polynomial time, provided the screening precedes any continuous relaxation. On the resulting graph, Hall's theorem characterizes when every driver can be served, and limits on the number of vehicles a site can host keep the problem exact, whereas limits on total site power do not. We further observe that sequential least-cost booking is a serial dictatorship and therefore, by a classical result, Pareto efficient for the drivers it serves. Any saving from joint assignment is thus a net gain that necessarily makes some served drivers worse off, and caps on individual losses can be imposed exactly by pruning pairs. On archived inventories of three American cities with simulated demand, joint assignment lowers the cost of serving the same drivers by 2.0 to 8.9 percent, while about two in five of them lose. With requests arriving over time, batching for about five minutes lowers the full cost per request by 1 to 2 percent once the confirmation delay is valued.
\end{abstract}


\section{Introduction}\label{sec:intro}
\IEEEPARstart{T}{he} problem of choosing a charging station is familiar to the driver of an electric vehicle (EV). A nearby station may charge a low price for electricity but have a long queue or a slow outlet. A faster station farther away may be worth the trip for a driver who needs a large charge before returning to work. Station recommendation systems weigh such preferences \cite{habbal2024preference} and can take account of changing requests and station conditions \cite{khan2026clustering}. A fleet dispatcher faces a further question, even with the same information about prices, power and waiting times. If two drivers can reach the same fast port and only one of them can also reach a slower one, which driver should receive the fast port? The decision can determine whether both vehicles get a charge at all.

We consider a platform that receives the requests together and reserves a particular port for each vehicle it can serve. Each driver specifies how much energy the vehicle needs, how far it can travel to charge and the maximum bill for that energy. The decision comes after tariffs are posted and before port reservations are confirmed. By comparing the costs of serving exactly the drivers who obtained a sequential booking, we can ask what coordination saves them while keeping a separate account of the requests that went unfilled. Holding requests to decide them jointly also delays confirmation, and when ports turn over during the day, the allocation changes how many vehicles each port can serve.

Charging navigation systems already combine traffic and power-system data to route drivers to fast-charging stations \cite{Guo2014navigation}, and smart-charging controllers schedule the power drawn by vehicles that are already plugged in \cite{Wang2016survey,Lee2021ACN}. Port reservation sits between the two: it fixes which vehicle plugs in where, and therefore how long each vehicle occupies a socket and how a site's energy is shared.

Where the network itself is short of sockets, the same graph of acceptable pairs locates the drivers who compete for too few ports and the sites where one more socket would help; a companion study develops that planning use and the power a shared site connection can deliver \cite{Companion2026}. Here we focus on the dispatcher.

This paper makes the following contributions.
\begin{itemize}
\item We formulate port reservation with individual reachability, connector and budget constraints as a mixed-integer program and show that, once unaffordable or unreachable vehicle--port pairs are screened out, it is exactly a rectangular assignment problem, solvable in polynomial time (Section~\ref{sec:model}). The order matters: screening must precede any continuous relaxation, and Appendix~\ref{app:solver} gives an instance where relaxing first returns a strictly lower value that no reservation attains. On the screened graph, Hall's theorem characterizes full service, and limits on the number of vehicles per site preserve integrality through a network-flow formulation, whereas limits on total site power do not (Section~\ref{sec:weighted}).
\item We show that sequential least-cost booking is a serial dictatorship and hence, by a classical result adapted here to classes of identical sockets, Pareto efficient for the drivers it serves. The saving from joint assignment is therefore a net gain that necessarily leaves some served drivers worse off, and caps on individual losses can be imposed exactly by pruning pairs (Section~\ref{sec:fairness}).
\item We simulate arrivals over time with ports that turn over, identify the batching window that balances coordination against confirmation delay and the load at which batching costs service, and show that an occupancy price on fast ports recovers part of that service (Section~\ref{sec:rolling}).
\item We compare five booking policies on three archived city networks with simulated demand, decompose the saving into travel, energy, waiting and charging components, measure who gains and who loses, and test the results against connector restrictions, street routing, inventory classification and the slowdown of fast charging near a full battery (Sections~\ref{sec:design} and~\ref{sec:city}).
\end{itemize}
Section~\ref{sec:related} reviews related work, and Section~\ref{sec:discussion} discusses implications and limitations.

\section{Related work}\label{sec:related}

Charging-station recommendation has mostly been studied from the perspective of an individual driver or a taxi fleet. Tian et al.\ recommend stations to electric taxis in real time to reduce waiting \cite{Tian2016}, Cao et al.\ select stations for vehicles on the move using their anticipated charging reservations to minimize trip duration \cite{Cao2018}, and Habbal and Alrifaie rank stations by user preferences over several attributes \cite{habbal2024preference}. Learning-based methods recommend stations with multi-agent reinforcement learning \cite{ZhangW2021}, plan charging with deep reinforcement learning to reduce total charging time \cite{ZhangC2021}, or combine clustering with dynamic requests and station status \cite{khan2026clustering}. Navigation schemes couple station choice with the power system and traffic, through real-time grid and traffic data \cite{Guo2014navigation} or a hierarchical game between stations and vehicles \cite{TanWang2017}, and distributed scheduling balances queues among highway service stations \cite{Gusrialdi2017}. Moghaddam et al.\ design a charging strategy for networks that offer slow, fast and battery-swapping options \cite{Moghaddam2018}, and Shen et al.\ review optimization models for charging operations and infrastructure planning \cite{Shen2019}.

Closer to our setting, fleet charging-location assignment has been formulated as an optimization model \cite{martinez2025fleet}, and Rigas et al.\ compare offline and online mechanisms that allocate vehicles to stations, accounting for reachability and willingness to pay \cite{Rigas2022}. The Hungarian algorithm has been used to assign vehicles to battery-swapping stations \cite{You2017} and to pair drivers for vehicle-to-vehicle charge sharing \cite{Bakr2023}. In ride-hailing, requests are pooled over short intervals and assigned jointly \cite{AlonsoMora2017,Yan2020}, and the interval itself can be optimized \cite{Qin2021}. Online bipartite matching studies the opposite extreme, irrevocable decisions on arrival \cite{KVV1990}. Sequential least-cost booking is a serial dictatorship, a mechanism known to give Pareto-efficient allocations \cite{Abdulkadiroglu1998}; we use this to explain who gains and who loses from coordination. For planning, Campa\~na et al.\ use assignment to size charging facilities under urban traffic flows \cite{Campana2021}, Rizopoulos and Gkiotsalitis link infrastructure investment to electric bus schedules and energy prices \cite{Rizopoulos2025}, and Parent et al.\ study station expansion in Montr\'eal with a network-flow model built from observed sessions and trips \cite{Parent2024}; our companion study uses matching witnesses for siting \cite{Companion2026}.

Table~\ref{tab:literature} summarizes these approaches. Most recommend or schedule stations for individual drivers or fleets with heuristic, game-theoretic or learning methods that do not certify optimality. We instead decide port-level reservations for a batch of drivers with individual travel, connector and budget constraints, solve the problem exactly, measure who gains and who loses, and test the batching decision over time.

\begin{table}[t]
\centering\footnotesize\setlength{\tabcolsep}{3pt}
\caption{Representative Approaches to Charging-Station Choice and Allocation}
\label{tab:literature}
\begin{tabular}{@{}l>{\raggedright\arraybackslash}p{3.3cm}>{\raggedright\arraybackslash}p{3.55cm}@{}}
\toprule
Ref. & Decision & Approach\\
\midrule
\cite{Tian2016} & Station recommendation for electric taxis & Real-time recommendation system\\
\cite{Cao2018} & Station selection on the move & Reservation-aware selection scheme\\
\cite{Guo2014navigation} & Fast-charging navigation & Power-system and traffic data\\
\cite{TanWang2017} & Navigation and station pricing & Hierarchical game\\
\cite{Gusrialdi2017} & Highway station scheduling & Distributed cooperative control\\
\cite{Moghaddam2018} & Choice among charging options & Network charging strategy\\
\cite{Rigas2022} & Vehicle-to-station allocation & Offline and online mechanism design\\
\cite{ZhangW2021} & Station recommendation & Multi-agent reinforcement learning\\
\cite{ZhangC2021} & Charging planning & Deep reinforcement learning\\
\cite{habbal2024preference} & Station recommendation & Multi-attribute user preferences\\
\cite{khan2026clustering} & Station selection & Clustering and deep Q-network\\
\cite{martinez2025fleet} & Fleet charging-location assignment & Optimization model\\
This paper & Port reservation for a batch & Exact assignment with screening, rolling-horizon test\\
\bottomrule
\end{tabular}
\end{table}

\section{Optimization model}\label{sec:model}

We model a \emph{snapshot} of an urban charging system. At the moment the platform makes its decision, \(N\) vehicles have requested a charge and \(M\) ports can be reserved. The distinction between a station and a port matters here. A station, or \emph{site}, is a location, whereas a \emph{port} is an outlet that can charge one vehicle. A station with six ports at the same address offers six places to plug in. Once a port is promised to one driver, the platform cannot offer it to another in the same batch, even if the first driver has yet to arrive. We first formulate the decision when a reservation is required for every request, then consider what the platform can do if that proves impossible.

\subsection{Decision variables and data}\label{sec:variables}

For every port \(i=1,\dots,M\) and vehicle \(j=1,\dots,N\), let \(x_{ij}=1\) if the platform assigns port \(i\) to vehicle \(j\), and let it be zero otherwise. The binary matrix \(X=[x_{ij}]\) records the entire decision. Reading down a column tells a driver where to go. Reading across a row tells the operator which vehicle, if any, will occupy that port.

A port has a posted electricity tariff \(p_i\) in dollars per kilowatt-hour (kWh), a power rating in kilowatts (kW) and an expected wait \(w_i\) in minutes (min) before a previous user releases it. The wait represents activity outside the batch and can be zero. A vehicle requests \(q_j>0\) kWh and its driver sets a session budget \(b_j\geq0\) in dollars. The requested energy is fixed before dispatch, so at port \(i\) the bill would be \(p_iq_j\). A driver requesting 40~kWh with a \$20 budget can therefore accept a tariff of up to \$0.50/kWh, however much travel or waiting the driver might save at a more expensive station.

Two further inputs describe the journey. The one-way distance \(d_{ij}\) is measured from the request origin to the port in kilometres (km). The driver's travel limit \(r_j\), also in kilometres, can reflect remaining battery range or a limit on how far the driver will go for a charge. We measure the whole trip to the port. Measuring extra distance relative to a planned journey would also require its destination or route. Let \(h_{ij}=1\) if the vehicle has a compatible connector and the driver has access to the port, and zero otherwise. The \emph{reachability matrix} \(R\) then has entry \(R_{ij}=1\) precisely when \(h_{ij}=1\) and \(d_{ij}\le r_j\). This records where a vehicle can physically go and plug in before we consider its driver's spending limit.

\subsection{Time out of service and the cost of an assignment}\label{sec:cost}

A fleet pays for the electricity it buys and also loses the use of a vehicle while it is at the station. For a compatible vehicle and port, let \(P_{ij}>0\) be their effective charging power, the lower of the port rating and the power the vehicle can accept. At a constant rate it takes \(q_j/P_{ij}\) hours to supply the requested energy. The expected \emph{time out of service at the station}, in minutes, is therefore
\begin{equation}
 T_{ij}=w_i+60\,q_j/P_{ij}.
 \label{eq:stationtime}
\end{equation}
This expression distinguishes the time until a socket becomes free from the time the vehicle actually draws power. We value both at \(v\) dollars per minute and the journey at \(a\) dollars per kilometre. The resulting \emph{generalized cost}, expressed in dollars for each assignment, is
\begin{equation}
\begin{aligned}
 c_{ij}&=a\,d_{ij}+p_i q_j+v\,T_{ij}\\
       &=\underbrace{a\,d_{ij}}_{\text{travel}}
       +\underbrace{p_i q_j}_{\text{energy bill}}
       +\underbrace{v\,w_i}_{\text{waiting}}
       +\underbrace{60v\,q_j/P_{ij}}_{\text{charging}} .
\end{aligned}
 \label{eq:cost}
\end{equation}
Only the electricity bill is charged against the session budget; travel and time enter the choice among offers that pass this test. The energy demand, tariff and effective power are known at assignment, so $c_{ij}$ is a coefficient of the objective. This treats the effective power as a property of the vehicle and the port alone, which holds when a site connection can supply every socket it occupies at full rating. A companion study relaxes that assumption and measures what a shared connection delivers~\cite{Companion2026}.

Scaling all weights by a positive constant preserves the ordering of assignments, so giving the bill a coefficient of one leaves only \(a\) and \(v\) to choose. In the city experiments \(a=\$2.52/\mathrm{km}\) and \(v=\$0.224/\mathrm{min}\), or \$13.44 per vehicle-hour, so at equal tariffs a port one kilometre farther away is worthwhile if it saves more than 11.25 minutes of waiting and charging. For perspective, the United States Department of Transportation values personal travel at \$20.10 per person-hour in 2024 dollars and the marginal operating cost of a light-duty vehicle at \$0.56 per mile \cite[Tables~A-2 and A-4]{USDOTtime2026}. Applied to a trip to the port and back at an average urban speed of about 22~km/h, these values give approximately our \$2.52 per kilometre of one-way distance. Appendix~\ref{app:calib} examines the choice of weights, and Section~\ref{sec:sensitivity} varies them.

\subsection{Objective function and constraints}\label{sec:objective}

When every request is to receive a port, the platform seeks a joint assignment of least total cost. The following binary linear program, a special case of a mixed-integer linear program (MILP), states this decision directly.
\begin{equation}
\begin{aligned}
 \min_X\quad &\sum_{j=1}^{N}\sum_{i=1}^{M}c_{ij}x_{ij} \\
 \text{subject to}\quad
 &\sum_{i=1}^{M}x_{ij}=1 && \forall j,\\
 &\sum_{j=1}^{N}x_{ij}\le 1 && \forall i,\\
 &x_{ij}\le R_{ij} && \forall i,j,\\
 &\sum_{i=1}^{M}p_i q_j x_{ij}\le b_j && \forall j,\\
 &x_{ij}\in\{0,1\} && \forall i,j.
\end{aligned}
\label{eq:fullprogram}
\end{equation}
The first two constraints give each vehicle one port and reserve each port for at most one vehicle. The reachability and budget constraints express the two conditions on the driver's choice. A vehicle must be able to get there and plug in, and the price of the energy it needs must be within the driver's budget. Even with more ports than vehicles, two drivers whose only acceptable choice is the same socket cannot both be served.

\subsection{From acceptable ports to a matching}\label{sec:matching}

Before committing any reservations, the dispatcher can determine which offers each driver could accept. Let \(\rho_j=b_j/q_j\) be driver \(j\)'s highest acceptable unit price. The physically available ports form the set \(\mathcal R_j=\{i:R_{ij}=1\}\), while the ports within budget form \(\mathcal B_j=\{i:p_i\le\rho_j\}\). We record their intersection, the set \(\mathcal A_j\) of acceptable ports, in column \(j\) of a binary matrix \(A\).
\begin{equation}
 \mathcal A_j=\mathcal R_j\cap\mathcal B_j,\qquad
 A_{ij}=\mathbf{1}\{i\in\mathcal A_j\}=R_{ij}\,\mathbf{1}\{p_i\le\rho_j\}.
 \label{eq:edges}
\end{equation}
Here \(\mathbf{1}\) is one when the condition in braces holds and zero otherwise.

\begin{proposition}[Screening individual reservations]\label{proposition:delta_equivalence}
If a vehicle can reserve at most one port and reservations are binary, its physical eligibility and budget constraints in~\eqref{eq:fullprogram} are equivalent to \(x_{ij}\le A_{ij}\) at every port \(i\).
\end{proposition}
\begin{proof}
If vehicle \(j\) is assigned to port \(i\), its bill is \(p_iq_j\). The chosen pair satisfies both constraints exactly when \(R_{ij}=1\) and \(p_iq_j\le b_j\), which by~\eqref{eq:edges} means \(A_{ij}=1\). All other entries in the column are zero. If the vehicle is unassigned, the column and the bill are both zero, and the budget is nonnegative.
\end{proof}

\begin{figure}[t]
\centering
\begin{minipage}[c]{.25\columnwidth}\centering\footnotesize\setlength{\arraycolsep}{2pt}
\(\begin{array}{c|cccc}
 i\backslash j&1&2&3&4\\ \hline
 1&1&0&0&0\\
 2&0&1&0&0\\
 3&1&1&0&1\\
 4&0&0&1&1\\
 5&0&0&0&1
\end{array}\)
\end{minipage}\hfill
\begin{minipage}[c]{.72\columnwidth}\centering
\includegraphics[width=\linewidth,trim=17bp 5bp 17bp 5bp,clip]{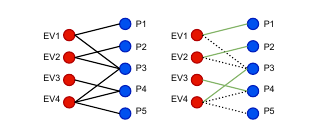}
\end{minipage}
\caption{From a matrix of acceptable pairs to a reservation. Rows index ports (P1--P5) and columns index vehicles (EV1--EV4). The middle panel draws every acceptable pair as an edge. The solid edges in the right panel show one assignment of all four vehicles, leaving port~5 free.}
\label{fig:acceptable}
\end{figure}

Each zero in \(A\) therefore lets us discard a decision variable, leaving one per acceptable pair. Fig.~\ref{fig:acceptable} illustrates the reduction with four vehicles and five ports. Each one in the matrix is an edge in a \emph{bipartite graph}, with vehicles on one side and ports on the other. A \emph{matching} selects edges that share neither a vehicle nor a port. Vehicle~3 can use only port~4, while vehicle~4 can use ports~3, 4 or 5. Reserving port~3 for vehicle~4 leaves port~4 for vehicle~3 and makes it possible to serve both.

\begin{theorem}[Full-service feasibility under individual charging budgets]\label{theorem:feasibility}
The charging model~\eqref{eq:fullprogram} admits a reservation for every vehicle if and only if the graph of acceptable pairs defined by~\eqref{eq:edges} has a matching that covers every vehicle. Equivalently, writing \(\Gamma(S)=\bigcup_{j\in S}\mathcal A_j\) for the set of ports acceptable to at least one vehicle in \(S\), full service is possible exactly when
\begin{equation}
 |\Gamma(S)|\ge |S|\qquad\text{for every group of vehicles }S.
 \label{eq:hall}
\end{equation}
\end{theorem}
\begin{proof}
By Proposition~\ref{proposition:delta_equivalence}, every pair selected in a feasible assignment is an edge of the graph. The assignment gives one edge to each vehicle and uses each port at most once, so the selected edges form a matching. Conversely, placing a one on each edge of a matching that covers the fleet and zeros elsewhere gives a matrix \(X\) satisfying~\eqref{eq:fullprogram}. Hall's theorem~\cite{Hall1935,LovaszPlummer1986} gives the second equivalence, since a matching covering the vehicle side exists exactly when~\eqref{eq:hall} holds.
\end{proof}

Taking \(S\) to be the whole fleet in Theorem~\ref{theorem:feasibility} gives \(N\le|\Gamma(S)|\le M\) whenever full service is feasible. When \(N=M\), such a reservation uses every port. A smaller group can still prevent full service, as when three vehicles can use only the same two sockets.

\subsection{Maximum service and the Hungarian algorithm for minimum cost}\label{sec:weighted}

The same graph describes the operator's choices when full service is impossible. Write \(E=\{(j,i):A_{ij}=1\}\) for its edges, and let \(K^\star\) be the greatest number of pairs in any matching. For a prescribed number \(K\) of reservations, the cost problem is
\begin{equation}
 \min_{\substack{F\subseteq E\ \text{a matching}\\ |F|=K}}\ 
       \sum_{(j,i)\in F} c_{ij}.
 \label{eq:assignment}
\end{equation}
We take \(K=K^\star\), giving priority to service and then to cost among assignments that serve that many drivers. When \(K^\star=N\) this is precisely~\eqref{eq:fullprogram}. Otherwise, two assignments with the same service count may leave different drivers without a port. We therefore record the number and energy requests of those drivers and compare booking costs on a common group in Section~\ref{sec:measure}.

The dispatcher can find \(K^\star\) without enumerating the groups in Hall's condition. The Hopcroft--Karp algorithm finds a maximum matching in a number of operations bounded by a constant times \((|E|+N+M)\sqrt{N+M}\), where \(|E|\) is the number of acceptable pairs \cite{Hopcroft1973}. If everyone can be served, the remaining decision is the classical linear assignment problem on those pairs. If \(K^\star<N\), add \(N-K^\star\) artificial ports that any driver may take at zero cost. A maximum matching can be completed with these ports, and every assignment covering all vehicles must use all of them because at most \(K^\star\) real ports can be assigned. The cheapest complete assignment therefore gives the cheapest way to serve \(K^\star\) drivers once the artificial reservations are discarded.

Kuhn's Hungarian algorithm solves the minimum-cost assignment problem exactly \cite{Kuhn1955,Burkard2009}. It maintains potentials \(u_j\) for vehicles and \(\beta_i\) for ports with \(u_j+\beta_i\le c_{ij}\) on every acceptable pair, whose sum bounds the cost of any complete assignment from below, and stops when a complete assignment attains equality on every chosen pair. Zero-cost artificial vehicles fill out the square when ports are in surplus, and modern implementations take \(O(n^3)\) operations for \(n\) rows \cite{Burkard2009}. The city experiments use a related method that works directly with rectangular matrices \cite{Jonker1987,Crouse2016} (Appendix~\ref{app:solver}). The model treats $c_{ij}$ as fixed at booking. If a site connection cannot supply all its occupied sockets at full rating, plugging in another vehicle lowers the power at a neighbouring one, the costs become interdependent and the problem leaves the assignment formulation~\cite{Companion2026}.

A limit on the number of vehicles that a site may host at once, for example the number of full-power sessions its connection can supply, keeps the problem exact. Routing each port's arc in the network-flow form of Appendix~\ref{app:solver} through a site node, whose arc to the sink carries that capacity, leaves a network matrix whose linear program has integral optima. A limit on total site energy or power is a knapsack constraint and loses this property.

Algorithm~\ref{alg:dispatch} summarizes one dispatch epoch. The static studies of Section~\ref{sec:city} apply it once to a batch; the rolling-horizon study of Section~\ref{sec:rolling} applies it every \(\Delta\) minutes to the requests still pending and the ports then free.

\begin{algorithm}[t]
\caption{One dispatch epoch: screening and maximum service at minimum cost}
\label{alg:dispatch}
\begin{algorithmic}[1]
\Require pending requests \(\mathcal J\) with \(q_j,\rho_j,r_j\); free ports \(\mathcal I\) with tariffs \(p_i\); distances \(d_{ij}\); compatibility \(h_{ij}\); costs \(c_{ij}\) from~\eqref{eq:cost}
\State Screen pairs: \(A_{ij}\gets h_{ij}\,\mathbf 1\{d_{ij}\le r_j\}\,\mathbf 1\{p_i\le\rho_j\}\), as in~\eqref{eq:edges}
\State \(K^\star\gets\) size of a maximum matching in the graph of \(A\) (Hopcroft--Karp)
\If{\(K^\star<|\mathcal J|\)} add \(|\mathcal J|-K^\star\) zero-cost artificial ports open to every request \EndIf
\State Solve the rectangular assignment with costs \(c_{ij}\) on acceptable pairs (Jonker--Volgenant)
\State Discard artificial reservations and confirm the others
\State Keep unserved requests pending until their patience expires
\end{algorithmic}
\end{algorithm}

\subsection{Example: allocating one fast port}\label{sec:example}
Suppose two drivers request 25 and 50~kWh. A 7.2~kW port charges \$0.35/kWh and a 50~kW port charges \$0.50/kWh. Both vehicles can draw the stated power at either port, and both offers meet their travel and spending limits. The smaller request originates 0.5~km from the fast port and 2~km from the slow one, with the distances reversed for the larger request. If both ports are free, assigning each driver to the nearer one gives 1~km of total travel, a \$30.00 combined bill and 446.7 vehicle-minutes of charging.

Giving the fast port to the larger request instead adds 3~km of travel and \$3.75 to the bill, while saving 178.3 vehicle-minutes of charging. At the reference distance weight the extra travel costs \$7.56, so the exchange is worthwhile when a vehicle-hour is worth more than \$3.81. At the city reference weights the two drivers save \$28.64 together. If the driver needing 50~kWh sets a price cap below \$0.50/kWh, the budget test removes this exchange before the costs are compared.

\section{Three-city study design}\label{sec:design}

\subsection{Charging supply and the geographic footprint}

We use an archived United States Department of Energy station inventory and three city files with geocoded addresses \cite{AFDCstations,census_geocoder,opencage}. The national file records access status and counts of Level~1 and Level~2 alternating-current (AC) ports and direct-current (DC) fast ports. AC equipment uses the vehicle's onboard charger, whereas a DC fast charger supplies direct current to the battery \cite{USA_charger_types}. 
The archived inventories accompany the analysis code in the replication package.

We retain the footprint of each city file, including one Alexandria record in Washington. Four Chicago records lack usable coordinates. We match the remaining records to the national extract by street address, city, state, five-digit postal code and port counts in each class. We keep only records at publicly accessible sites, and a public record enters the inventory at most as many times as it occurs in the source. In Chicago these two filters remove 13 ports at private sites and a further 114 public ports listed more than once, leaving 598 ports. Table~\ref{tab:inventory} gives the resulting supply, and Section~\ref{sec:access} tests how that choice affects dispatch.

\begin{figure*}[!t]
\centering
\includegraphics[width=\textwidth]{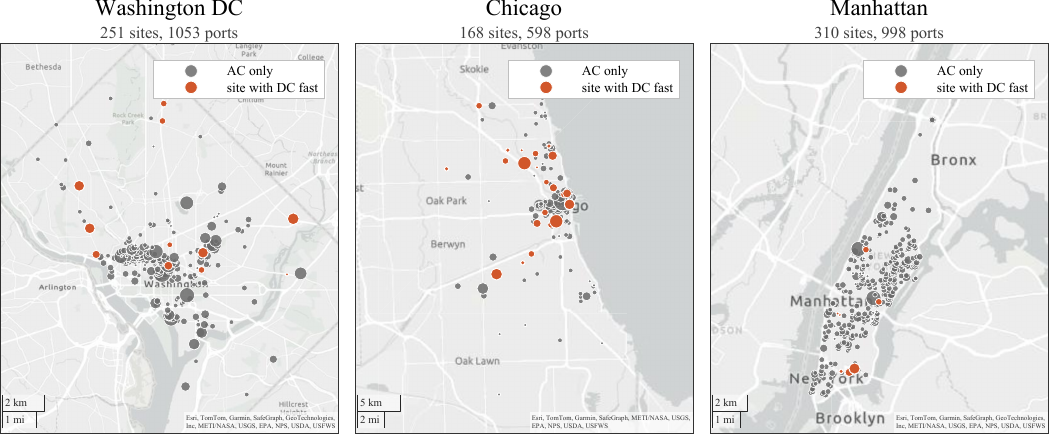}
\caption{The three charging networks, after the inventory is restricted to
publicly accessible ports and repeated public entries are capped at their
source counts. Marker area is proportional to the number of ports a site
hosts; orange marks a site with at least one DC fast port. Manhattan packs
310 sites into a narrow strip but holds only 20 fast ports, while Chicago
spreads 168 sites over a far wider area and hosts 102. Basemap: the \texttt{streets-light} map hosted by Esri, rendered with
MATLAB's \texttt{geobasemap}~\cite{geobasemap}; map data \copyright~Esri,
TomTom, Garmin, SafeGraph, GeoTechnologies Inc., METI/NASA, USGS, EPA, NPS,
USDA, USFWS.}
\label{fig:networks}
\end{figure*}

\begin{table}[t]
\centering\footnotesize\setlength{\tabcolsep}{3pt}
\caption{Public Charging Supply After Matching the Geocoded Records to the National Inventory}
\label{tab:inventory}
\begin{tabular}{lrrrrrrr}
\toprule
City & Records & Sites & L1 & L2 & DC & Ports & Requests\\
\midrule
Washington & 334 & 251 & 2 & 988 & 63 & 1,053 & 527 \\
Chicago & 230 & 168 & 0 & 496 & 102 & 598 & 299 \\
Manhattan & 340 & 310 & 0 & 978 & 20 & 998 & 499 \\
\bottomrule
\end{tabular}

\smallskip
\parbox{\columnwidth}{\footnotesize A site is a distinct coordinate. L1 and L2 count Level~1 and Level~2 alternating-current ports, DC counts direct-current fast ports, and Requests is the half-load batch, rounded with halves rounded up.}
\end{table}

Sockets at one coordinate share a posted tariff and an estimated wait but retain their individual equipment classes. A coordinate can represent colocated equipment or reflect the resolution of the address match. Fig.~\ref{fig:networks} maps the three resulting networks.

For Level~1 AC, Level~2 AC and DC fast ports we assign reference ratings of 1.4, 7.2 and 50~kW, respectively, within the corresponding equipment classes \cite{USA_charger_types}. Each vehicle's AC acceptance limit is equally likely to be 7.2 or 11~kW and its DC limit equally likely to be 50 or 100~kW. The lower of the vehicle limit and the port rating sets the charging rate. A 40~kWh request therefore takes about 333 minutes at 7.2~kW and 48 minutes at 50~kW. The reference fleet can connect to every equipment class. Section~\ref{sec:connectors} varies access to fast ports using the connector labels on the archived records.

\subsection{Requests, prices and external demand}

A \emph{layout} is one simulated set of vehicle origins, requests, tariffs and waits. We locate each origin relative to the mean coordinates of the retained records, drawing independent eastward and northward offsets from normal distributions with standard deviation 2.5~km. We retain origins inside the convex hull of the sites, the smallest convex region containing them all. This concentrates demand near the centre of the inventory. The city files contain no observed vehicle origins or session records from which to estimate this distribution.

In the main comparison each driver can travel up to 5~km in a straight line to reach a port. We draw energy demand \(q_j\) uniformly between 20 and 60~kWh and the price cap \(\rho_j\) uniformly between \$0.40 and \$0.60/kWh, independently. Setting \(b_j=\rho_jq_j\) means that a larger request takes longer to charge and costs more at a given tariff, while its probability of passing the budget test stays the same.

Each site receives a posted price drawn uniformly between \$0.303 and \$0.530/kWh, a range motivated by published United States charging-price data \cite{USA_charge_price}, and an external wait drawn uniformly between zero and 15 minutes. A dispatcher may choose a site with a shorter wait, but no policy can reschedule those earlier customers. Prices and waits are independent of one another and of power in the reference case. Section~\ref{sec:sensitivity} examines the effect of linking them.

\subsection{Booking policies}\label{sec:policies}

The first four policies use two criteria for choosing a port and two ways of committing a reservation. A \emph{batch} policy sees all pending requests. A \emph{sequential} policy sees the current driver, the site conditions and which ports are free, then confirms a reservation before seeing the next request. Every policy respects the same eligibility rules.

\emph{Batch cost} solves \eqref{eq:assignment} at maximum service using \eqref{eq:cost}. \emph{Batch distance} also serves as many vehicles as possible, but among those assignments chooses one with the least total travel. Different sockets at one coordinate have exactly the same distance. Having chosen the site for each served driver, we assign sockets within each site to minimize generalized cost while preserving all the site choices. A distance-based operator therefore receives the benefit of allocating faster equipment to the requests for which it matters most within a site.

\emph{Sequential cost} gives the current driver the free, eligible port of least generalized cost. \emph{Nearest} chooses the nearest site with an available eligible port. The fifth policy, \emph{cheapest of three}, first identifies the three nearest distinct eligible sites with a free port, or all such sites if fewer remain. From these it selects the site with the lowest tariff. Both site-based policies use the least-cost free socket at their chosen site. A common random ordering of ports breaks remaining equal port scores, and the sequential methods process requests in the same random order. If every eligible port is occupied, the request is recorded as unserved.

\subsection{Performance measures}\label{sec:measure}

For each layout and arrival order, let \(S_b\) be the drivers served by policy \(b\), \(K_b=|S_b|\) their number and \(C_b\) their combined generalized cost. The service rate is \(K_b/N\). The service deficit \((K^\star-K_b)/N\) counts the additional reservations that were possible as a fraction of all requests. We also record the total energy requested by the drivers left without a reservation.

To compare costs for the same people, we give the batch optimizer all the original ports and ask it to serve exactly \(S_b\). Write \(C^\star(S_b)\) for the minimum cost of serving that group, also called the \emph{served cohort}. The percentage saving relative to the competing policy is
\begin{equation}
 G_b=100\,\frac{C_b-C^\star(S_b)}{C_b}.
 \label{eq:saving}
\end{equation}
Because the policy's own assignment is feasible for \(S_b\), \(G_b\) is nonnegative for a nonempty group. It measures what the dispatcher could have saved the same drivers by allocating their ports jointly.

To account for the drivers left without a reservation, we assign a hypothetical fallback cost \(h\) to each unsuccessful referral. The full-fleet cost per requested vehicle is \(L_b(h)=[C_b+h(N-K_b)]/N\). We consider \(h=\$50,\$100,\$200\). One can interpret \(h\) as the total generalized consequence of seeking an alternative arrangement, including lost time, if no port is offered within this network. We have no field estimate of that consequence, so we vary it after dispatch with all assignments held fixed.

City policy comparisons use 30 independent layouts per city, load and operating scenario. We average three common random processing orders within each layout. For a paired 95\% bootstrap interval, we draw 30 entire layouts with replacement, keep all policies and orders together, and take the 2.5th and 97.5th percentiles of 2,000 replicate estimates \cite{Efron1979}. These intervals capture simulation variability conditional on the archived inventory and the demand model.

\section{Three-city results}\label{sec:city}

\subsection{Policy comparison at half load}

We begin with one request per two installed ports, to the nearest whole request. A seemingly generous aggregate supply can conceal quite different choices at the driver's starting point. Under the same origin-generation rule, the mean straight-line distance to the nearest installed site, before checking prices or reservations, is 0.44~km in Washington, 0.74~km in Chicago and 0.24~km in Manhattan. Applying both the 5~km radius and the budget cap leaves a driver an average of 157, 73 and 199 distinct eligible sites, respectively, before other drivers take their ports. These are averages over simulated requests and layouts. They describe how the sampled demand meets the local supply, rather than the distribution of trip origins for all city residents.

\begin{table*}[t]
\centering\small\setlength{\tabcolsep}{3.6pt}
\caption{Half-load dispatch on 30 layouts per city. The percentage saving uses the same served drivers as the policy in that row, with a paired 95\% bootstrap interval. Travel, energy bill, waiting time and charging time are per driver served by the indicated policy.}
\label{tab:city}
\begin{tabular}{llrr rrrr}
\toprule
City & Policy & Served & Same-driver saving & Travel & Bill & Wait & Charge\\
    &        & (\%) & (\%) & (km) & (\$) & (min) & (min)\\
\midrule
Washington & Batch cost & 100.00 & --- & 1.02 & 14.70 & 6.5 & 283.3 \\
 & Batch distance & 100.00 & 7.68 [7.48, 7.88] & 0.72 & 16.18 & 7.5 & 309.5 \\
 & Sequential cost & 99.94 & 4.63 [4.55, 4.71] & 1.30 & 14.64 & 6.7 & 297.9 \\
 & Nearest & 99.99 & 8.21 [8.02, 8.41] & 0.82 & 16.09 & 7.5 & 311.1 \\
 & Cheapest of three & 99.96 & 8.56 [8.16, 8.95] & 1.00 & 15.09 & 7.6 & 315.0 \\
Chicago & Batch cost & 99.94 & --- & 1.64 & 15.62 & 6.6 & 224.0 \\
 & Batch distance & 99.94 & 13.27 [12.96, 13.58] & 1.24 & 16.50 & 7.3 & 272.6 \\
 & Sequential cost & 99.48 & 8.92 [8.73, 9.12] & 1.94 & 15.40 & 6.7 & 252.0 \\
 & Nearest & 99.12 & 15.08 [14.78, 15.39] & 1.38 & 16.41 & 7.3 & 277.9 \\
 & Cheapest of three & 99.09 & 15.44 [14.96, 15.95] & 1.63 & 15.63 & 7.4 & 280.0 \\
Manhattan & Batch cost & 100.00 & --- & 0.50 & 14.64 & 6.5 & 316.6 \\
 & Batch distance & 100.00 & 3.99 [3.83, 4.15] & 0.31 & 16.23 & 7.5 & 327.1 \\
 & Sequential cost & 100.00 & 2.03 [1.98, 2.08] & 0.72 & 14.72 & 6.7 & 321.6 \\
 & Nearest & 100.00 & 4.03 [3.86, 4.18] & 0.34 & 16.13 & 7.5 & 327.2 \\
 & Cheapest of three & 100.00 & 3.42 [3.30, 3.55] & 0.48 & 15.14 & 7.5 & 327.5 \\
\bottomrule
\end{tabular}
\end{table*}

Table~\ref{tab:city} shows that both batch policies serve every driver in all 30 Washington and Manhattan layouts. In Chicago they serve everybody on 25 of 30 layouts and miss a single request on each of the other five. Four of those five contain an isolated vehicle, meaning one with no eligible port before any reservations have been made. In the fifth every driver has an option, yet a group competes for too few sockets. Across all Chicago layouts, joint assignment serves 99.94\% of requests, compared with 99.48\% for sequential cost, a paired difference of 0.46 percentage points [0.36, 0.56].

In Washington, sequential cost pays \$14.64 for energy per served driver, slightly less than the \$14.70 paid by batch cost, so a comparison of receipts would favour it. Including trips and plugged-in time reverses the ranking: sequential cost assigns trips of 1.30~km and occupies ports for 297.9 minutes per served driver, against 1.02~km and 283.3 minutes under batch cost. The tariff says little about how many vehicle-hours a particular driver saves on a fast socket.

Chicago shows the largest cost opportunity. Batch distance gives trips of 1.24~km on average, against 1.64~km under batch cost, but mean charging durations of 272.6 rather than 224.0 minutes, and reassigning its own drivers with the full objective saves 13.27\% [12.96, 13.58]. In Manhattan, even a rule that checks the prices of three nearby sites serves all sampled requests, although joint reassignment of its drivers can still save 3.42\%.

\subsection{Decomposition of the saving}\label{sec:components}

For the informed sequential rule the same-driver savings are \WashingtonSaving\% in Washington, \ChicagoSaving\% in Chicago and \ManhattanSaving\% in Manhattan. The order follows fast-port scarcity. Chicago has 102 fast ports at 22 sites, against 63 at 11 sites in Washington and 20 at seven sites in Manhattan, so at most 34\%, 12\% and 4\% of the respective batches could use one; Fig.~\ref{fig:networks} shows where those sites are. In Manhattan nearly every driver has a fast port within reach and budget, but only a few can reserve one in the same batch.

Table~\ref{tab:components} subtracts the four cost components of each sequential assignment and of its minimum-cost reassignment on the same drivers. A positive entry is a saving, and the small negative bill entries mean that the operator lowers the overall cost by buying slightly more expensive electricity. In Chicago, \$6.37 of the \$6.98 saving comes from shorter charging sessions, about 28 minutes per served driver at the reference valuation. The externally imposed wait contributes only a few cents.

\begin{table}[t]
\centering\footnotesize\setlength{\tabcolsep}{2.6pt}
\caption{Generalized-Cost Savings When Reassigning Exactly the Sequential Policy's Successful Requests}
\label{tab:components}
\begin{tabular}{lrrrrr rr}
\toprule
&\multicolumn{5}{c}{Class-specific power (\$/driver)}&\multicolumn{2}{c}{Common power}\\
\cmidrule(lr){2-6}\cmidrule(lr){7-8}
City & Travel & Bill & Wait & Charge & Total & \$/driver & \%\\
\midrule
Washington & 0.72 & -0.06 & 0.05 & 3.29 & 3.99 & 0.65 & 0.70 \\
Chicago & 0.81 & -0.22 & 0.03 & 6.37 & 6.98 & 0.79 & 0.82 \\
Manhattan & 0.56 & 0.07 & 0.05 & 1.14 & 1.83 & 0.68 & 0.74 \\
\bottomrule
\end{tabular}

\smallskip
\parbox{\columnwidth}{\footnotesize The last two columns repeat the entire policy comparison with every port rated at 7.2~kW. Components and totals are rounded independently.}
\end{table}

When every port supplies 7.2~kW, a driver's charging duration is the same wherever the vehicle plugs in. Repeating both policies on the same inputs now leaves only distance, tariff and external wait to distinguish the offers. Savings for the drivers served by sequential booking fall to \$0.65 in Washington, \$0.79 in Chicago and \$0.68 in Manhattan, below 1\% in every city (Table~\ref{tab:components}).

\subsection{Who gains and who loses}\label{sec:fairness}

Averages hide who pays for coordination. Under sequential cost booking each driver, in processing order, takes the cheapest free eligible port, which is a serial dictatorship \cite{Abdulkadiroglu1998}. Sockets with the same site and equipment class have the same eligibility and cost for every driver, so we group them into classes. Serial dictatorships are known to yield Pareto-efficient allocations \cite{Abdulkadiroglu1998}; the following proposition adapts this classical result to classes of identical sockets, and we include its short proof for completeness.

\begin{proposition}[Sequential booking is Pareto efficient]\label{prop:pareto}
Suppose that, for every driver, distinct socket classes have distinct generalized costs. Then no reassignment of the drivers served by sequential cost booking lowers one driver's cost without raising another's.
\end{proposition}
\begin{proof}
Let \(X\) be the sequential assignment and \(Y\) a reassignment of the same drivers in which no driver pays more. Suppose some driver holds different classes in \(X\) and \(Y\), and let \(k\) be the first such driver in processing order, with class \(x\) in \(X\) and \(y\) in \(Y\). Drivers processed before \(k\) hold the same classes in both assignments, and under \(Y\) they leave at least one socket of class \(y\) for \(k\). That socket was therefore free when \(k\) chose under \(X\), and it is eligible for \(k\). Since \(k\) chose the cheapest free eligible class and the costs of distinct classes differ, \(c_k(x)<c_k(y)\), so \(k\) pays more under \(Y\), a contradiction. Hence every driver keeps its class and every cost is unchanged.
\end{proof}

\begin{figure*}[!t]
\centering
\includegraphics[width=\textwidth]{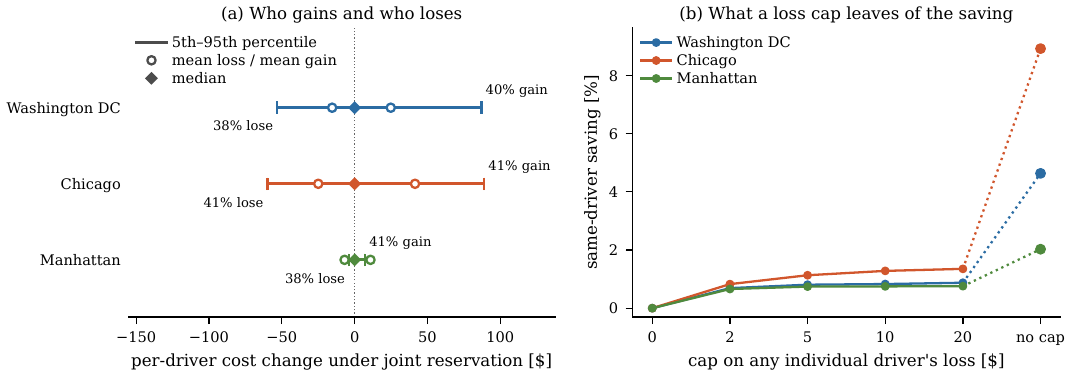}
\caption{Per-driver outcomes of joint reservation, against sequential
least-cost booking on the drivers it serves. (a) The spread of the cost
change, positive when the driver pays less; the bar covers the 5th to the
95th percentile, the open circles are the mean loss among those who lose and
the mean gain among those who gain. (b) The same-driver saving that survives
when no driver may lose more than the amount on the horizontal axis. At a
zero cap the saving is exactly zero, as Proposition~\ref{prop:pareto}
requires; most of it only appears once drivers are allowed to lose tens of
dollars.}
\label{fig:outcomes}
\end{figure*}

With continuous distances and tariffs the condition holds with probability one. Joint reassignment therefore cannot help one served driver without hurting another, and Table~\ref{tab:fairness} and Fig.~\ref{fig:outcomes}(a) show the split. In every city about 40\% of the sequentially served drivers gain and 38--41\% lose. The median change is zero in all three: between a sixth and a quarter of the drivers keep an equally costly reservation, and the rest divide almost evenly between gaining and losing, so the redistribution falls entirely on the tails. Those tails differ by an order of magnitude across the networks. In Washington the fifth and ninety-fifth percentiles of the change lie near $-\$53$ and $\$87$, and in Chicago near $-\$60$ and $\$89$, whereas in Manhattan they are $-\$4$ and $\$7$, because the network holds only twenty fast ports and almost no driver can be moved onto one. The drivers who gain request about 45~kWh on average and those who lose about 36~kWh, because the reassignment moves large requests onto fast ports and smaller ones off them. Total gains exceed total losses by a factor of 1.7 in each city, so compensating every displaced driver from the saving is affordable. A cap on individual losses can also be imposed exactly, by removing the pairs that would exceed it before solving the assignment. Fig.~\ref{fig:outcomes}(b) traces the whole curve. With a zero cap the saving is exactly zero in every layout, as Proposition~2 predicts; with a \$5 cap it falls to 0.81\%, 1.13\% and 0.74\%, and raising the cap to \$20 adds almost nothing. The curve is almost flat between \$2 and \$20 and rises steeply only once the cap is lifted altogether, so the saving does not accumulate gradually as drivers are allowed to lose more. Most of the benefit of coordination thus rests on moving some drivers off fast ports at a cost of tens of dollars each, which a platform would need to compensate or price.

\begin{table}[t]
\centering\footnotesize\setlength{\tabcolsep}{2.1pt}
\caption{Per-Driver Outcomes When the Sequentially Served Drivers Are Reassigned Jointly at Half Load}
\label{tab:fairness}
\begin{tabular}{@{}lrrrrrrrr@{}}
\toprule
& \multicolumn{2}{c}{Drivers (\%)} & \multicolumn{2}{c}{Mean (\$)} & Lose & Gain/ & \multicolumn{2}{c}{Capped (\%)}\\
\cmidrule(lr){2-3}\cmidrule(lr){4-5}\cmidrule(l){8-9}
City & gain & lose & gain & loss & $>\$5$ & loss & \$5 & \$20\\
\midrule
Washington & 39.7 & 37.9 & 24.79 & 15.45 & 10.8 & 1.68 & 0.81 & 0.87 \\
Chicago & 41.4 & 40.7 & 41.47 & 24.98 & 21.4 & 1.69 & 1.13 & 1.35 \\
Manhattan & 41.2 & 38.5 & 10.89 & 6.91 & 4.1 & 1.69 & 0.74 & 0.76 \\
\bottomrule
\end{tabular}

\smallskip
\parbox{\columnwidth}{\footnotesize Differences below half a cent count as ties. Mean gain and loss are over the drivers who gain or lose; \emph{Lose \(>\$5\)} is the percentage of drivers who lose more than \$5. \emph{Capped} is the same-driver saving when no driver may lose more than \$5 or \$20; uncapped savings are in Table~\ref{tab:city}. Averages over 30 layouts and three processing orders.}
\end{table}

\subsection{Connector restrictions and the cost of delay}\label{sec:connectors}

A station record identifies connector types, but it does not identify the vehicles or adapters that will turn up on a given day. We therefore repeat the half-load dispatches with two hypothetical fleets. In the first, each driver independently has a one-half chance of being unable to use fast ports on records that list Tesla as their only connector. In the second, this restriction applies to every driver. Other ports, including those on mixed-connector records, keep their original eligibility. Prices, requests, locations and arrival orders are paired across cases. Within each case we compare costs on exactly the drivers who obtained a sequential reservation.

\begin{table*}[t]
\centering\small\setlength{\tabcolsep}{4pt}
\caption{Cost savings and the value of booking time at half load. The fast-port count gives the private ones out of all fast ports. Gains are the percentage reduction in cost for the sequentially served drivers in each fleet. The difference is in percentage points, with a paired 95\% interval. Delay is the additional mean time before confirmation that would use up the reference fleet's saving.}
\label{tab:connectors}
\begin{tabular}{lrrrrrr}
\toprule
 & Fast & \multicolumn{3}{c}{Same-driver gain (\%)} & Restricted minus & Break-even\\
\cmidrule(lr){3-5}
City & ports & Reference & Mixed & Restricted & reference (points) & delay (min)\\
\midrule
Washington & 36/63 & 4.63 & 4.43 & 2.56 & $-2.08$ [$-2.15, -2.00$] & 17.8 \\
Chicago & 56/102 & 8.92 & 8.48 & 4.98 & $-3.95$ [$-4.09, -3.80$] & 31.2 \\
Manhattan & 14/20 & 2.03 & 2.00 & 1.14 & $-0.88$ [$-0.92, -0.85$] & 8.2 \\
\bottomrule
\end{tabular}
\end{table*}

Chicago retains the largest saving even with all drivers restricted (Table~\ref{tab:connectors}), and in the mixed fleet its saving of 8.48\% stays close to the reference 8.92\%, because a compatible driver can still use a scarce fast port when another vehicle cannot.

A dispatcher who holds requests for joint booking introduces a wait before confirming the port. If each extra minute costs the same \(v=\$0.224\) as a minute at the station, a mean saving of \(g\) dollars per served driver pays for \(g/v\) additional minutes: 17.8 in Washington, 31.2 in Chicago and 8.2 in Manhattan. The relevant delay is the average extra wait relative to immediate confirmation, which depends on when requests arrive, so the threshold alone cannot fix a booking interval; Section~\ref{sec:rolling} tests this directly.

\subsection{Service under increasing load}\label{sec:load}

\begin{figure*}[t]
\centering
\includegraphics[width=.88\textwidth]{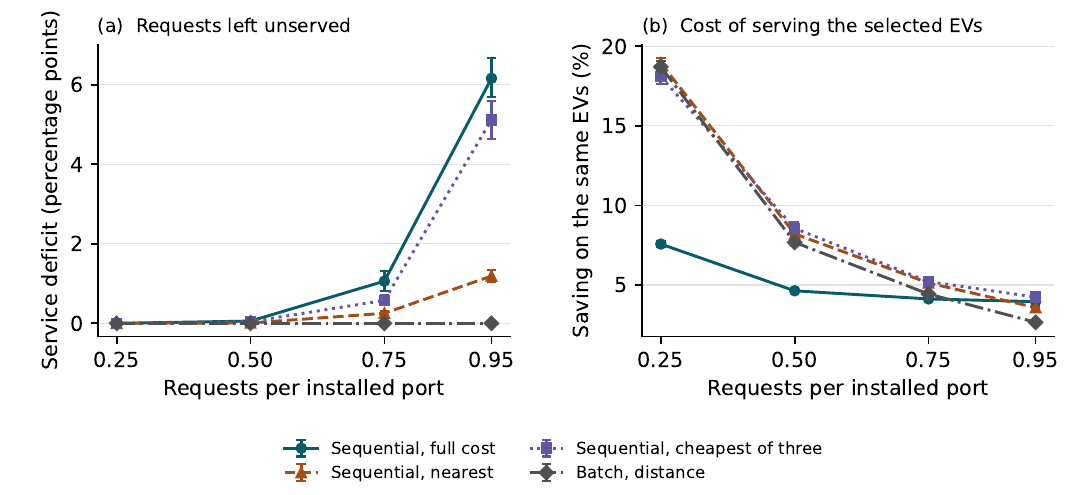}
\caption{Service and cost as the Washington batch grows, over 30 layouts per setting. Panel (a) gives the service deficit in percentage points relative to the maximum matching. Panel (b) gives the percentage saving available on the same drivers served by each policy. Error bars show paired 95\% bootstrap intervals.}
\label{fig:load}
\end{figure*}

Holding Washington supply fixed at 1,053 ports, we increase the batch from 263 to 527, 790 and 1,000 requests. Fig.~\ref{fig:load} compares the fraction of requests a policy misses despite sufficient eligible ports with the savings still possible for its successful drivers.

With 263 requests, both batch rules and the sequential full-cost rule serve every driver in every sampled batch. Reassignment can still lower the cost of serving those drivers by 7.57\% [7.40, 7.74]. With 1,000 requests, both batch rules again find room for everyone, while sequential cost serves \HighGreedyService\%. Its deficit of \HighGap{} percentage points [\HighGapLo, \HighGapHi] corresponds to about 62 missed vehicles and \HighUnservedEnergy{} megawatt-hours (MWh) of requested energy per batch. On the drivers who do receive an offer, the remaining opportunity to save through reassignment is \HighSaving\%. Thus the smaller percentage saving at high load coexists with a larger opportunity to offer service. As more ports fill, there are fewer interchangeable choices for the served drivers, and committing the wrong one can prevent a later driver from receiving any offer at all.

At the highest load, sequential cost exceeds batch cost by \$\FallbackFifty, \$\FallbackHundred{} or \$\FallbackTwoHundred{} per \emph{requested} vehicle when a missed referral is valued at \$50, \$100 or \$200, each additional \$100 adding \$6.16 per requested vehicle to the difference.

\subsection{Forecast-guided sequential booking}\label{sec:forecast}

An operator who must confirm reservations immediately might still anticipate competition. In a deliberately favourable test, the operator solves an assignment linear program on one independent forecast batch, with sockets of the same site and class forming capacity groups, and adds each group's shadow price \cite[Sec.~5.6]{Boyd2004} to the dispatch score of its free sockets. At half load this rule leaves a same-driver saving of 1.99\% [1.94, 2.05], against 4.63\% for the ordinary sequential rule; at high load it raises service only from \HighGreedyService\% to 94.09\%. A forecast thus recovers much of the cost saving but little of the service gain.

\subsection{Dispatch over time and the batching window}\label{sec:rolling}

To test batching when requests arrive over time and ports turn over, we simulate six hours of arrivals on the Washington network: a Poisson stream of 1.5, 3.0 or 4.5 requests per minute (540, 1,080 or 1,620 per run) with the demand model, eligibility and generalized cost of Section~\ref{sec:design}. A reservation holds the port from confirmation until the vehicle leaves, covering travel at 22~km/h, the site's external wait and the charging session; a port can be reserved only when free, and a request unassigned 30 minutes after arrival is lost and incurs the \$100 fallback. Immediate booking gives each arrival the least-cost free eligible port and retries unassigned requests every minute. Batched booking applies Algorithm~\ref{alg:dispatch} every \(\Delta\) minutes, for \(\Delta\) between 1 and 30. Each served driver's cost adds the confirmation delay at \(v=\$0.224\)/min, and we draw 30 layouts per rate.

\begin{figure*}[t]
\centering
\includegraphics[width=.88\textwidth]{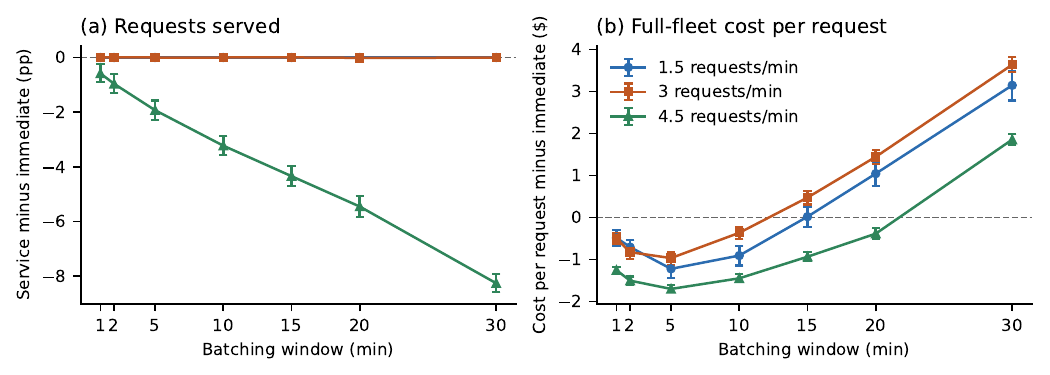}
\caption{Rolling-horizon dispatch on the Washington network over six hours of arrivals, relative to immediate booking. Panel (a) gives the change in the percentage of requests served; the curves for 1.5 and 3.0 requests per minute coincide at zero. Panel (b) gives the change in full-fleet cost per request, including the confirmation delay valued at \$0.224/min and a \$100 fallback for each lost request. Error bars show paired 95\% bootstrap intervals over 30 layouts.}
\label{fig:rolling}
\end{figure*}

\begin{table*}[t]
\centering\small
\caption{Rolling-Horizon Dispatch on the Washington Network over Six Hours of Arrivals}
\label{tab:rolling}
\begin{tabular}{lrrrrrr}
\toprule
Arrivals & Batching & Served & Mean delay & Cost per served & Cost per & Change vs immediate\\
(per min) & window (min) & (\%) & (min) & driver (\$) & request (\$) & (\$) [95\% CI]\\
\midrule
1.5 & Immediate & 100.00 & 0.0 & 56.98 & 56.98 & --- \\
 & 1 & 100.00 & 0.5 & 56.39 & 56.50 & $-$0.48 [$-$0.67, $-$0.30] \\
 & 5 & 100.00 & 2.5 & 55.20 & 55.76 & $-$1.22 [$-$1.44, $-$1.02] \\
 & 10 & 100.00 & 5.0 & 54.95 & 56.08 & $-$0.90 [$-$1.14, $-$0.68] \\
 & 20 & 100.00 & 10.0 & 55.79 & 58.03 & $+$1.05 [$+$0.76, $+$1.32] \\
 & 30 & 100.00 & 15.0 & 56.77 & 60.13 & $+$3.15 [$+$2.79, $+$3.48] \\
3 & Immediate & 99.99 & 0.0 & 74.30 & 74.31 & --- \\
 & 1 & 99.99 & 0.5 & 73.70 & 73.81 & $-$0.50 [$-$0.63, $-$0.37] \\
 & 5 & 99.99 & 2.5 & 72.78 & 73.34 & $-$0.96 [$-$1.13, $-$0.83] \\
 & 10 & 99.99 & 5.0 & 72.82 & 73.95 & $-$0.36 [$-$0.52, $-$0.22] \\
 & 20 & 99.98 & 10.0 & 73.51 & 75.75 & $+$1.45 [$+$1.28, $+$1.61] \\
 & 30 & 99.98 & 15.0 & 74.58 & 77.95 & $+$3.64 [$+$3.46, $+$3.82] \\
4.5 & Immediate & 97.95 & 2.3 & 81.68 & 82.57 & --- \\
 & 1 & 97.37 & 1.5 & 80.48 & 81.32 & $-$1.25 [$-$1.33, $-$1.17] \\
 & 5 & 96.02 & 3.3 & 79.33 & 80.87 & $-$1.70 [$-$1.78, $-$1.61] \\
 & 10 & 94.72 & 5.7 & 78.80 & 81.12 & $-$1.45 [$-$1.54, $-$1.35] \\
 & 20 & 92.50 & 10.3 & 78.43 & 82.19 & $-$0.38 [$-$0.50, $-$0.26] \\
 & 30 & 89.69 & 15.0 & 79.27 & 84.43 & $+$1.86 [$+$1.73, $+$1.98] \\
\bottomrule
\end{tabular}

\smallskip
\parbox{\textwidth}{\footnotesize Cost per served driver excludes the confirmation delay. Cost per request adds the delay valued at \$0.224/min for served drivers and a \$100 fallback for each lost request. Means and paired bootstrap intervals over 30 layouts.}
\end{table*}

Fig.~\ref{fig:rolling} and Table~\ref{tab:rolling} show the trade-off. Immediate booking costs \$56.98, \$74.31 and \$82.57 per request at the three rates, and the mean confirmation delay under batching at the two lower rates is about half the window. At 1.5 and 3.0 requests per minute essentially every request is served. A five-minute window gives the lowest full cost per request, \$1.22 [1.02, 1.44] and \$0.96 [0.83, 1.13] below immediate booking, or about 2.1\% and 1.3\%. Beyond 10 to 15 minutes the delay costs more than coordination saves. The dynamic saving is smaller than the static one because each fast port now serves several vehicles over six hours, and immediate booking reuses it as soon as it frees. At 4.5 requests per minute the network saturates. Immediate booking serves 97.95\% of requests, and batching serves fewer: 0.58 percentage points fewer with one-minute windows and 1.93 with five-minute windows. Batching gives fast ports to larger requests, 46.7 against 39.9~kWh on average at five minutes, and their longer sessions let each fast port serve fewer vehicles. Batching still lowers the full cost per request for windows up to 20 minutes, by \$1.70 [1.61, 1.78] at five minutes, because the cheaper and faster sessions of the served drivers outweigh the fallback cost of the extra lost requests. An operator that values service above cost near saturation should therefore batch briefly. The static finding that coordination adds reservations (Section~\ref{sec:load}) does not carry over directly once ports turn over.

Because the loss comes from how fast ports are shared, we add to the dispatcher's scores, but not to the reported costs, an occupancy price of \(\kappa\) dollars per minute that a fast port is held. Tuned on ten separate layouts, \(\kappa=\$0.8\)/min maximizes service at 4.5 requests per minute; above about \$0.9/min a fast port costs more than the charging time it saves and drivers abandon it. On the evaluation layouts the price narrows the service loss with five-minute windows from 1.93 to 1.29 percentage points [0.94, 1.64], and with one-minute windows from 0.58 to an insignificant 0.21 [\(-0.10\), 0.50], while lowering the cost slightly. At 1.5 requests per minute, however, the same price raises the cost by \$1.07 per request, because fast ports are not scarce there. Like a congestion charge, the price should respond to scarcity.

\begin{table}[t]
\centering\footnotesize\setlength{\tabcolsep}{2.6pt}
\caption{Robustness of Five-Minute Batching at 4.5 Requests per Minute}
\label{tab:rollrobust}
\begin{tabular}{@{}lrrr@{}}
\toprule
& Change in & Change in cost per & Mean delay\\
Case & served (pp) & request (\$) [95\% CI] & (min)\\
\midrule
Base & $-1.93$ & $-1.70$ [$-1.78$, $-1.61$] & 2.3 / 3.3\\
Three-hour warm-up & $-1.25$ & $-6.01$ [$-6.17$, $-5.85$] & 8.6 / 4.2\\
Patience 15 min & $-1.64$ & $-1.75$ [$-1.84$, $-1.64$] & 1.0 / 2.8\\
Patience 60 min & $-1.57$ & $-1.61$ [$-1.73$, $-1.48$] & 4.0 / 4.5\\
Fallback \$50 & $-1.93$ & $-2.66$ [$-2.86$, $-2.45$] & 2.3 / 3.3\\
Fallback \$200 & $-1.93$ & $+0.23$ [$-0.12$, $0.59$] & 2.3 / 3.3\\
Price \$0.8/min & $-1.29$ & $-1.76$ [$-1.84$, $-1.67$] & 2.3 / 3.3\\
\bottomrule
\end{tabular}

\smallskip
\parbox{\columnwidth}{\footnotesize Changes are five-minute batching minus immediate booking, paired over 30 layouts. Immediate booking serves 97.95\% of requests in the base case, 88.87\% after the warm-up and 96.29\% and 99.77\% with patience of 15 and 60 minutes. Mean delay is the confirmation delay of served drivers under immediate and batched booking. The warm-up case simulates nine hours and counts requests arriving after the third.}
\end{table}

Table~\ref{tab:rollrobust} varies the remaining assumptions at 4.5 requests per minute; at the two lower rates every variant leaves service unchanged and the five-minute saving between \$0.96 and \$1.37 per request. A three-hour warm-up within a nine-hour horizon removes the initial filling of an empty network. In that steady state immediate booking serves only 88.9\% of requests and a queue forms: retried requests wait 8.6 minutes on average and, served first come first served, often take ports that others needed. Batching serves 1.25 points fewer but assigns the backlog jointly at maximum service, halves the mean wait and lowers the full cost per request by \$6.01, the static matching argument reappearing over time. Batching windows must also stay shorter than drivers' patience: with 15-minute patience, 20-minute windows lose about a quarter of all requests. Finally, the verdict at saturation depends on how a lost request is valued; with a \$200 fallback, five-minute batching only breaks even.

\subsection{Sensitivity to prices, geography and charging power}\label{sec:sensitivity}

Table~\ref{tab:sensitivity} varies one assumption at a time on the Washington half-load layouts. The circuity case multiplies straight-line distances by 1.30; the two-centre case draws demand around two sites at opposite quartiles of the network's main axis, with a standard deviation of 1.5~km.

The equipment case allows Level~2 ratings of 7.2, 11 and 19~kW and DC ratings of 50, 100 and 150~kW, equally probable within each class, subject to the vehicle's acceptance limit. The linked-price case gives sites with fast equipment a price premium and makes cheaper sites more likely to have a longer external wait, as specified in the replication package. We also test a lower valuation of station time (\$0.10/min) and a higher valuation of distance (\$5.04/km). Finally, the taper case lets the average power of a DC session fall linearly with the request, from the nominal effective power at 20~kWh to half of it at 60~kWh, a simple stand-in for the slowdown of fast charging at high states of charge.

\begin{table*}[t]
\centering\small\setlength{\tabcolsep}{4pt}
\caption{Washington at half load under alternative assumptions. Complete counts layouts where all requests could be served. The sequential full-cost rule's deficit is in percentage points of all requests. Saving is the percentage available by reassigning its served drivers. Brackets give paired 95\% bootstrap intervals. Bill and charge describe the sequential rule.}
\label{tab:sensitivity}
\begin{tabular}{lcrrrr}
\toprule
Case & Complete & Service deficit & Same-driver saving & Bill (\$) & Charge (min)\\
\midrule
Base & 30/30 & 0.06 [0.03, 0.08] & 4.63 [4.55, 4.71] & 14.64 & 297.9 \\
Circuity & 30/30 & 0.24 [0.18, 0.31] & 4.54 [4.43, 4.66] & 14.71 & 298.3 \\
Two centres & 30/30 & 0.00 [0.00, 0.01] & 3.77 [3.70, 3.85] & 14.77 & 300.4 \\
Variable power & 30/30 & 0.04 [0.02, 0.07] & 6.95 [6.84, 7.05] & 14.64 & 246.1 \\
DC charging taper & 30/30 & 0.06 [0.03, 0.08] & 3.50 [3.42, 3.57] & 14.64 & 300.5 \\
Linked tariffs and waits & 30/30 & 0.00 [0.00, 0.01] & 4.24 [4.16, 4.30] & 14.35 & 297.9 \\
Lower time value & 30/30 & 0.06 [0.03, 0.09] & 4.57 [4.48, 4.65] & 14.59 & 297.9 \\
Higher distance value & 30/30 & 0.04 [0.02, 0.07] & 5.03 [4.94, 5.11] & 14.86 & 298.0 \\
\bottomrule
\end{tabular}
\end{table*}

The same-driver saving ranges from 3.50\%, when DC sessions taper, to 6.95\%, with varied power within classes. The taper removes about a quarter of the base saving but leaves the service deficit unchanged. Multiplying distances by 1.30 raises the sequential service deficit from 0.06 to 0.24 percentage points.

\subsection{Street-network distances}\label{sec:roads}

We repeat the half-load comparison on the District Department of Transportation's Roadway Block network \cite{DDOTblocks}: 24,772 directed links between 8,638 endpoints, with origins and stations joined to their nearest endpoint by straight access legs (median 61 and 51~m). Distance is the two access legs plus the shortest directed path, and a pair without such a path is ineligible. Maximum service falls from 100\% to 98.67\%, a paired loss of 1.33 percentage points [1.18, 1.46], entirely because three to eleven drivers per layout become isolated, and sequential service falls from 99.94\% to 98.37\%. The same-driver saving is 4.58\% [4.47, 4.68], close to the straight-line estimate, and 4.55\% [4.44, 4.66] when access is restricted to the largest \emph{strongly connected component} of 8,505 mutually reachable endpoints, where maximum service recovers to 99.97\%. The cost opportunity is thus robust to the street model, whereas who can be served depends on where vehicles can enter the network.

\subsection{Inventory classification}\label{sec:access}

Inventory curation moves the estimates by about one point, and nearly all of it comes from duplication rather than from access. On the all-listed Chicago footprint the saving is 8.81\% when every one of its 725 geocoded port entries is treated as accessible, 8.69\% when only publicly accessible sites are kept, and 7.77\% when repeated public entries are also limited to their source counts, a paired change of $-1.05$ points $[-1.16, -0.94]$ across the two steps. These arms hold the demand of the all-listed footprint fixed, so they are not directly comparable with Table III, which is computed on the curated inventory at its own half load. In Manhattan, removing 63 private ports lowers the saving from 2.99\% to 2.03\%. All Washington ports matched public records, and the ranking of the cities is unchanged.

\section{Discussion and conclusions}\label{sec:discussion}

A platform with requests in hand can lower the generalized cost of the drivers it serves by changing who gets a fast socket: by 2.0--8.9\% in a single batch, mostly through shorter charging sessions, and by about 1--2\% per request over time with a five-minute window. The saving falls when fewer vehicles can use fast connectors or fast sessions taper, almost disappears when every port has the same power, and is lost to the confirmation delay for windows beyond 10 to 20 minutes.

The gain is a net one. Sequential booking is Pareto efficient for the drivers it serves, and joint assignment leaves about two in five of them worse off, by \$7 to \$25 on average depending on the city. Gains exceed losses by a factor of about 1.7, so compensation is affordable, but caps on individual losses remove most of the saving. A platform adopting joint booking would need a transfer or pricing rule, for instance a discount for drivers moved away from a fast port, together with a model of drivers who decline the recommendation.

Scarcity changes what the dispatcher needs from the model. Once range or connector restrictions differ, the shortage can lie among drivers whose only option is the same socket; in one batch of 1,000 Washington requests, maximum matching recovers about 62 reservations that immediate least-cost booking misses. When ports turn over, batching near saturation serves 1--2 points fewer requests than immediate booking because it gives fast ports to large requests; an occupancy price on fast ports, applied only under scarcity, recovers part of that loss. In a congested steady state, by contrast, joint assignment of the backlog halves the wait and lowers the full cost per request by about 7\%. Field tests of these effects require records of declined as well as successful requests.

For planning, the same graph locates the drivers who compete for too few sockets and the sites where one more would help, and a companion study shows that the value of such a socket also depends on the power the site connection can share among its vehicles \cite{Companion2026}.

Several simplifications bound these conclusions. Demand is simulated: origins, energy requests, budgets, tariffs, external waits and the arrival stream follow stated distributions rather than observed sessions, and the intervals describe simulation variability under those assumptions. Drivers are assumed to accept the recommended port, and cancellations and no-shows are ignored. The occupancy price is a single tuned constant rather than a price that adapts to congestion. Charging power is constant within a session except in the taper case, and the power available at a site does not depend on how many vehicles share it. Natural extensions are validation on operator records of requested, declined and completed sessions, compensation and acceptance models, and site power limits inside the assignment.

\section*{Acknowledgment}
The station inventories, Washington street data, analysis code, numerical results and reproduction instructions are openly available in the project repository \cite{ChargingArchive2026}. During the preparation of this work the authors used Claude (Anthropic) and Codex (OpenAI) to improve the language and presentation of the manuscript and to assist with the implementation of the analysis code in Sections~\ref{sec:design} and~\ref{sec:city}. The authors reviewed and verified all content, including every derivation and numerical result, and take full responsibility for the article.

\appendices
\section{Calibration of the travel and time weights}\label{app:calib}

To examine the weights we generate 15 feasible synthetic networks with 500 ports and 250 vehicles each. Distances are independent and uniform on [0,10]~km, each pair is eligible with probability 0.6, waits are uniform on [0,15]~min, tariffs on [\$0.303, \$0.530]/kWh and requests on [20,60]~kWh, and powers of 1.4, 7.2 and 50~kW occur in proportion to the Washington port counts. On a logarithmic \(13\times13\) grid with \(a\in[0.25,4]\)~\$/km and \(v\in[0.05,1]\)~\$/min we record the mean distance, bill and station time per driver. Every tested setting is Pareto efficient among those tested \cite{Marler2010}. Normalizing each outcome to [0,1] over the grid and minimizing the Euclidean distance to the ideal point selects \(a=3.175\) and \(v=0.287\), with a mean distance of 0.195~km, bill of \$14.70 and station time of 288.6~min, against grid minima of 0.134~km, \$14.40 and 286.7~min. The city reference weights \(a=2.52\) and \(v=0.224\) are the grid points one step below on each axis and imply nearly the same exchange between distance and station time (11.25 versus 11.06 minutes per kilometre). Section~\ref{sec:sensitivity} shows that lowering \(v\) to \$0.10/min or doubling \(a\) changes the Washington saving by less than half a percentage point.

\section{Exact solution and computational checks}\label{app:solver}

\subsection{Screening before relaxation}

Proposition~\ref{proposition:delta_equivalence} uses the fact that an integer reservation assigns a vehicle to at most one port. If we instead allow fractional assignments to several ports, the budget row in~\eqref{eq:fullprogram} limits their average bill. The order of screening and relaxation therefore matters, as a one-vehicle example illustrates.

A driver who needs 1~kWh has a \$2 budget and can reach ports charging \$1 and \$3 per kWh, with respective generalized costs of \$11 and \$3. Only the first is within budget, giving an optimal cost of \$11. If the solver keeps both ports and allows a half assignment to each, the averaged bill is \$2 and the objective falls to \$7. Screening out the \$3 tariff first restores the correct \$11 optimum even with continuous variables.

After screening, represent~\eqref{eq:assignment} as a network with a source connected to each vehicle, one arc per acceptable vehicle--port pair, and an arc from every port to a sink. Each arc carries at most one unit. Only the vehicle--port arcs carry costs \(c_{ij}\). Sending \(K\) units from source to sink represents \(K\) reservations for integer flows. The node--arc incidence matrix is \emph{totally unimodular}, meaning all square submatrices have determinant zero, one or minus one. With integer capacities and \(K\), the resulting linear program (LP) has an integral optimum whenever feasible \cite{ahuja1993network,Schrijver2003}. Fractional optima can still occur as averages of tied assignments.

\subsection{Solvers and checks}

For the city comparisons we use \texttt{linear\_sum\_assignment} in SciPy~1.7.3 \cite{scipy2020}, a modified Jonker--Volgenant algorithm for rectangular assignment \cite{Jonker1987,Crouse2016}, with infinite costs on ineligible pairs and, when full service is impossible, the \(N-K^\star\) zero-cost artificial ports of Section~\ref{sec:weighted}. The Washington investment comparisons instead give each driver a private artificial port costing more than \(N\) times the largest real cost, so that serving one more driver always outweighs any change in real costs, while assignments with equal service are ranked by real cost alone. Every run checks eligibility, distinct ports and service counts and reconstructs the cost components in~\eqref{eq:cost}.

An independent MATLAB R2020a implementation with Optimization Toolbox 8.5 gave the same minimum for the original MILP, the screened LP and the assignment solver on 12 feasible instances with individual budgets, and the dense Hungarian implementation and the filtered integer program agreed within \(10^{-6}\) dollars on all 1,494 of 1,500 random instances that admitted full service. In 30 timing instances with up to 200 vehicles and 240 ports, median end-to-end times were below 10~ms for the assignment solvers and below 161~ms for the integer program on an eight-thread Intel Mac. The studies used Python 3.9.7, SciPy 1.7.3 and Matplotlib 3.4.3, and all trials are archived in the replication package \cite{ChargingArchive2026}.

\bibliographystyle{IEEEtran}
\bibliography{charging_refs}

\end{document}